\documentclass[11pt]{amsart}
\usepackage{multirow}
\usepackage{graphicx}
\usepackage{amsmath}
\usepackage[mathscr]{eucal}
\usepackage{dsfont}
\usepackage{xypic}
\usepackage{graphicx}
\usepackage{url}
\usepackage{hyperref}

\newcommand{\repl}[3]{\ensuremath{#1[#2\to#3]}}

\newcommand{\M}{\ensuremath{\mathscr{M}}}

\theoremstyle{plain}
\newtheorem{theorem}{Theorem}[section]
\newtheorem{proposition}[theorem]{Proposition}
\newtheorem{conjecture}[theorem]{Conjecture}
\newtheorem{lemma}[theorem]{Lemma}
\newtheorem{corollary}[theorem]{Corollary}

\theoremstyle{definition}
\newtheorem{definition}[theorem]{Definition}
\newtheorem{remark}[theorem]{Remark}
\newtheorem{convention}[theorem]{Convention}
\newtheorem{notation}[theorem]{Notation}

\theoremstyle{remark}

\newcommand{\cocoa}{{\hbox{\rm C\kern-.13em o\kern-.07em C\kern-.13em
      o\kern-.15em A}}}

\begin{document}

\title[Maximum independent sets of de Bruijn graphs]{The
  maximum independent sets of de Bruijn graphs of diameter 3}

\author[D. A. Cartwright]{Dustin A. Cartwright}

\author[M. A. Cueto]{Mar\'ia Ang\'elica Cueto}
\address{Department of Mathematics\\
  University of California, Berkeley\\
  970 Evans Hall \#3840\\
  Berkeley, CA 94720-3840}
\email{\{dustin,macueto\}@math.berkeley.edu}

\author[E. A. Tobis]{Enrique A. Tobis $^\dag$}
\address{Department of Ophthalmology, Children's Hospital Boston,
  Center for Brain Science and Swartz Center for Theoretical
  Neuroscience, Harvard University, 300 Longwood Avenue, Boston,
  Massachusetts 02115}
\email{Enrique.Tobis@childrens.harvard.edu}

\begin{abstract}
  The nodes of the de Bruijn graph $B(d,3)$ consist of all strings of
  length $3$, taken from an alphabet of size $d$, with edges between
  words which are distinct substrings of a word of length~$4$.  We
  give an inductive characterization of the maximum independent sets
  of the de Bruijn graphs $B(d,3)$ and for the de Bruijn graph of
  diameter three with loops removed, for arbitrary alphabet size.  We
  derive a recurrence relation and an exponential generating function
  for their number. This recurrence allows us to construct
  exponentially many comma-free codes of length~3 with maximal
  cardinality.
\end{abstract}
\thanks{$^\dag$Corresponding author.\\D.~A. Cartwright was supported by
  the National Science Foundation grants DMS-0354321, DMS-0456960, and
  DMS-0757236.  M.~A. Cueto was supported by a UC Berkeley
  Chancellor's Fellowship.  E.~A. Tobis was supported by a CONICET
  doctoral fellowship, CONICET PIP 5617, ANPCyT PICT 20569 and UBACyT
  X042 and X064 grants.}
\keywords{de Bruijn graphs, maximum independent sets}
\subjclass[2000]{05C69,05A15}
\maketitle

\section{Introduction}
\label{sec:intro}

For any positive integers $d$ and~$D$, the \emph{de Bruijn} graph
$B(d,D)$ is the directed graph whose $d^D$~nodes consist of all the
$D$-digit words from the alphabet $\{0,\ldots,d-1\}$.  There is a directed edge
 from a word $x = x_1\ldots x_D$ to $y = y_1\ldots y_D$ if and
only if $x_2\ldots x_D = y_1\ldots y_{D-1}$. These graphs were
introduced in~\cite{debruijn}, under the name of \emph{$T$-nets}.
Since then, de Bruijn graphs have been used in several contexts,
notably as a network
topology~\cite{bermond,ptp,mukherjee}, 
and for building protein-binding microarrays~\cite{bulyk}.

\begin{figure}[h]
  \centering
  \includegraphics[scale=.8]{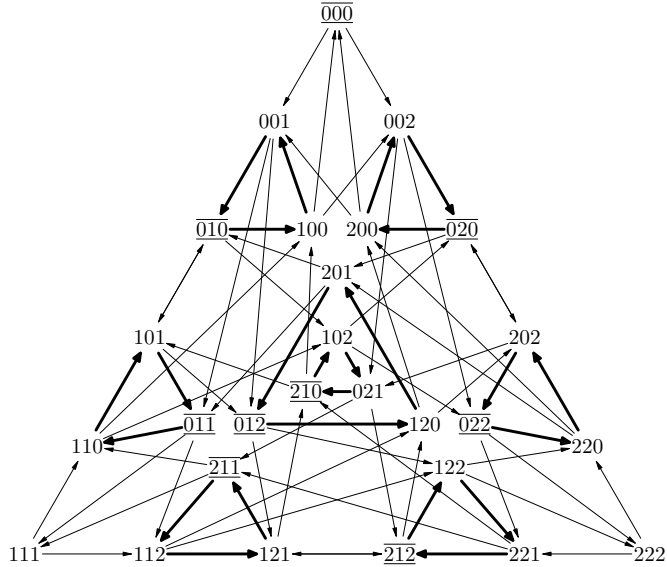}
  \caption{The de Bruijn graph $B(3,3)$, with the loops on $000$,
    $111$ and $222$ removed. The highlighted nodes belong to one of
    the 42 possible MISs in $B(3,3)$.  Bold arrows indicate edges
    under the shift function $\theta$ defined in~\eqref{eq:1}.
  }
  \label{fig:example}
\end{figure}

We concern ourselves with the maximum independent sets of these
graphs, previously studied in \cite{shibata,lichiardopol}. The graph
$B(d,D)$ contains $d$ nodes of the form $x\ldots x$, which have an
edge to themselves. In a slight abuse of notation, we will refer to
such a node as \emph{the loop~$x$}. Notice that a loop cannot be in
any independent set of $B(d,D)$, and therefore we call these sets
\emph{loop-less maximum independent sets} (LMISs).  The maximum
independent sets of the subgraph of $B(d,D)$ obtained by removing the
edges 
$x \ldots x\to
x\ldots x$ are called \emph{maximum independent sets} (MISs).
Figure~\ref{fig:example} depicts $B(3,3)$ with an MIS highlighted.

A natural question to ask is what is the \emph{stable size} of $B(d,D)$ for
arbitrary $d$ and $D$, i.e. the sizes of an MIS and a loop-less MIS. This
questions was studied in~\cite{lichiardopol}. 
Lichiardopol defined $\alpha(d,D)$ to be the size of an MIS with loops
and $\alpha^*(d,D)$ to be the size of a loop-less
MIS~\cite{lichiardopol}.
For $D$ a prime at least $3$, he proved the inequalities
\begin{equation}\label{eq:LichBound}
  \alpha(d,D) \leq \frac{(D-1)(d^D - d)}{2D} + 1 \quad \text{ and } \quad
  \alpha^*(d,D) \leq \frac{(D-1)(d^D - d)}{2D}.
\end{equation}
He then showed that in fact, equality holds for $D$ equal to~$3$, $5$
or~$7$ and conjectured that the same is true for all odd
primes~$D$. More precisely, it suffices to show that the conjecture
holds for $d=2$ and fixed $D$ to prove it for arbitrary $d$. As a byproduct of his work, we conclude that any MIS of $B(d,D)$ has at most two loops. 

In the case of $D=3$, we give a complete recursive characterization of
the maximum independent sets of $B(d,3)$.
To do so, we give four functions which extend an
MIS in $B(d,3)$ to an MIS in $B(d+1,3)$ or $B(d+2,3)$ (Definitions~\ref{def:f}, \ref{def:f'},
\ref{def:g} and \ref{def:g'}). Our main result
is that every maximum independent set in $B(d,3)$ can be formed by
beginning with an MIS in $B(1,3)$ or $B(2,3)$ and successively
applying our four functions and permuting the alphabet.
Moreover, since the sequence of functions and permutations is unique
up to certain transpositions, we can compute the number of
MISs, which corresponds to the Sloane sequence \texttt{A052608}~\cite{sloane}:

\vspace{1ex}
\noindent\textbf{Theorem \ref{theo:main}.} If we let $a_d$ be the number of
maximum independent sets of $B(d,3)$, then $a_d$ has exponential
generating function
\begin{equation*}
  \sum_{d=1}^\infty \frac{a_d t^d}{d!} = \frac{t+t^2}{1-2t-t^2}.
\end{equation*}
\vspace{1ex}

\noindent In addition, we prove that the number of loop-less
maximum independent sets has the same generating function (Theorem~\ref{theo:loop-less}).

A loop-less maximum independent set in $B(d,3)$ is a maximum comma-free code of
length~$3$. Comma-free codes were introduced by Crick, Griffith, and Orgel as a
hypothetical encoding of amino acid sequences in DNA~\cite{crick}, and further
generalized in~\cite{golomb,jiggs,golomb2}. A \emph{comma-free code} is a
set~$S$ of $D$-digit words such that if $x_1\ldots x_D$ and $y_1 \ldots y_D$ are
in~$S$, then no substring of $x_2 \ldots x_D y_1 \ldots y_{D-1}$ is in~$S$.
In~\cite{golomb} and~\cite{eastman}, it was shown that a comma-free code of
length $D = 3$ could have as many as $(d^3 - d) / 3$ elements by giving the same
example (up to permuting the alphabet), namely the code consisting of all words
$x_1 x_2 x_3$ such that $x_1 < x_2 \geq x_3$.  In contrast to this single
example, our results give an explicit construction of exponentially many
equivalence classes of maximum comma-free codes (Theorem~\ref{theo:comma-free}).

For $D=2$ and $d \geq 4$, the maximum independent sets have size
$\alpha(d,2) = \alpha^*(d,2) = \lfloor
d^2/4\rfloor$~\cite[Prop.~5.1]{lichiardopol}, and the same analysis as
in that proof shows that number of maximum independent sets of
$B(d,2)$ is $\binom{d}{d/2}$ if $d$ is even and $2\binom{d}{(d-1)/2}$
if $d$ is odd.

On the other hand, for $D > 3$, even small values of~$d$ yield de
Bruijn graphs with a large number of maximum independent sets in $B(d,
D)$.  For example, using the computer algebra system
\cocoa~\cite{cocoa}, we found out that there are $1$ and $44$ maximum
independents sets of $B(1,5)$ and $B(2,5)$, respectively. However, we
know that there are at least $210492$ maximum independent sets of
$B(3,5)$. 
This rapid growth means that the maximum independent sets in $B(3,5)$
cannot be produced from smaller independent sets using only
permutations of the alphabet and a handful of functions.  

We conjecture that an analogue of
Theorem~\ref{theo:main}(\ref{item:unique-sequence}) for $D > 3$ would
require starting with MISs in $B(d, D)$ for all $d <D$. Moreover, as
it occurs for diameter three, we would also need functions taking a
maximum independent set in $B(d, D)$ to one in $B(d + k, D)$ for all
$k < D$. Finding explicit formulas for these functions would require
knowledge of the sets $B(k,D)$ for $k<D$. To summarize:


\begin{conjecture}
  Let $D$ be a fixed odd prime number.  Then the exponential generating
  function of the number of maximum independent sets of $B(d, D)$ is
  the ratio of two polynomials, each of degree $D-1$. 
\end{conjecture}

The rest of this paper is organized as follows.  In
Section~\ref{sec:second-results}, we define two functions $f$ and~$f'$
that take a maximum independent set of $B(d,3)$ to a maximum
independent sets of $B(d+1,3)$. Likewise, we construct another two
functions $g$ and~$g'$ that take a maximum independent set of $B(d,3)$
to a maximum independent set of $B(d+2,3)$.  In
Section~\ref{sec:symmetries-bd-3} we compute the stabilizers of the
maximum independent sets produced by $f$, $f'$, $g$ and $g'$ under the
action of the symmetric group $\mathds{S}_d$, and show that the
functions take disjoint orbits to disjoint orbits.  In
Section~\ref{sec:maintheorem}, we prove our main theorems. In
Section~\ref{sec:loopless-miss}, we give a bijection between the
maximum independent sets and the loop-less maximum independent sets of
$B(d,3)$, from which we conclude that their numbers coincide.



\section{Inductive Construction of Maximum Independent Sets}
\label{sec:second-results}

In this section, we present two pairs of combinatorial operations that
transform a maximum independent set in the de Bruijn graph $B(d,3)$
into a maximum independent set in either $B(d+1,3)$ or
$B(d+2,3)$. 

\begin{convention}
  Throughout this paper, we let $\mathbf{[d]}$ stand for the set
  $\{0,\ldots,d-1\}$.
\end{convention}

Essential to the structure of the de Bruijn graph $B(d,3)$ are the
cycles under the shift function~$\theta$, defined as
\begin{equation}
  \label{eq:1}
  \theta\colon V(B(d,3)) \to V(B(d,3)) \qquad \theta(x y z) = y z x,
\end{equation}
where $V(B(d,3))$ denotes the set of nodes of the graph $B(d,3)$~\cite{lichiardopol}.
Note that the fixed points of~$\theta$ are exactly the loops
of~$B(d,3)$. On the other hand, if $x y z$ is not a loop, then $x y
z$, $\theta(x y z)$, and~$\theta^2(x y z)$ form a directed
$3$-cycle. 
In Figure~\ref{fig:example}, the
$\theta$-cycles are indicated by bold edges.

\begin{convention}
  Whenever we speak of \emph{cycles}, we mean the cycles induced by $\theta$.
\end{convention}

The action of $\theta$ induces a decomposition of the nodes
of~$B(d,3)$ into $(d^3 - d)/3$ cycles of length~$3$, and $d$ cycles of
length~$1$ (i.e.\ the loops). Each of these disjoint cycles
contributes at most one node to any independent set of
$B(d,3)$.

The following proposition explains the role played by a loop in a
maximum independent set of $B(d,D)$ and it shows that such a set can
have at most two loops. 

\begin{proposition}
  \label{prop:lessthanthreeloops}
  Let $D$ be an odd prime number, and let $S$ be a maximum independent
  set of $B(d,D)$ achieving the maximum possible size
  $\frac{(D-1)(d^D - d)}{2D} + 1$. 
  Then $S$ contains one or two loops. Moreover, if $a$ is a loop
  in~$S$ and $x$ 
  is any digit which is not a loop of~$S$, then the node
  $(ax)^{\frac{D-1}{2}}a$ is in~$S$.  If $S$ has two loops $a, b$,
  then, possibly after swapping $a$ and~$b$, $(ab)^{\frac{D-1}{2}}a$
  is in~$S$. Moreover, each cycle contributes exactly $(D-1)/2$ nodes
  to~$S$, except for one of the form $b^{D-2i}(ab)^{i}$, for some $1
  \leq i \leq (D-1)/2$, which only contributes $(D-1)/2 - 1$ nodes.
\end{proposition}
\begin{proof} The proof is contained in the proof of~\cite[Proposition
  4.3]{lichiardopol}.
\end{proof}
From the previous result we see that the loops of an MIS play a
special role. More precisely, if $S$ is a maximum independent set then
all the cycles of $B(d,D)$ of length $D$ contribute at most $(D-1)/2$
elements to $S$. If $D$ is an odd prime, and $S$ has only one loop,
then Lichiardopol's conjecture says that equality holds~\cite{lichiardopol}. If, on the
other hand, $S$ contains two loops $a, b$, then (up to swapping $a$
and~$b$) we
can assume that $(ab)^{\frac{D-1}{2}}a \in S$. Hence, all cycles of
$B(d,D)$ of length $D$ contribute at most $(D-1)/2$ nodes, except for
one cycle with more $b$'s than $a$'s, which contributes at most
$(D-1)/2 - 1$. Again, Lichiardopol's conjecture states that these
maximal contributions are achieved~\cite{lichiardopol}.

Since the conjecture holds for $D=3$, every cycle (with the possible
exception of the cycle of $bab$) contributes one element to any
maximum independent set. This motivates the following definition,
which will play an essential role in our inductive construction of
maximum independent sets of $B(d,3)$.



\begin{definition}\label{def:M_x(A)}
  Let $A$ be a set of nodes from $B(d,3)$. Let $x$ and $y$ be two
  digits in~$\mathbf{[d]}$. We say that \emph{$y$ appears between $x$ in
    $A$} if the node $x y x$ belongs to~$A$.
  We define $\M_x(A)$ as the set of digits which \emph{do not appear} between
  $x$ in $A$.
  We define $m_x(A)$ as the number of digits which do not appear between
  $x$ in $A$, i.e.\ $m_x(A) = \lvert\M_x(A)\rvert$.
\end{definition}

\begin{notation}\label{not:replaceAndLoops}
  If $w$ is a node in $B(d,3)$, we will denote by \repl{w}{x}{y} the node that
  results from replacing every occurrence of the digit $x$ by the
  digit $y$ in $w$.  We write $x \in w$ to mean that $x$ is one of the
  digits that appear in $w$.  

We denote by $L(S)$ the set of loops of a maximum independent set
$S$. We denote by $a$ the element of $L(S)$ such
  that $m_a(S) = 0$. We will refer to it as \emph{the
    distinguished loop}. 
If $S$ has another loop we denote it $b$.  This distinction will be
  extremely important for the construction of our four 
  operations on $B(d,3)$.
\end{notation}
We now define our first operation, sending a maximum independent set of
$B(d,3)$ to a subset of $B(d+1,3)$. Proposition~\ref{prop:fismis} will
show that this subset is a maximum independent set.
\begin{definition}
  \label{def:f}
  Let $S$ be a maximum independent set of $B(d,3)$. Following
  Notation~\ref{not:replaceAndLoops}, we define $f(S)
  \subset B(d+1,3)$ as the set $S\cup \bigcup_{i=1}^5 U_i(S)$, where
  \begin{align*}
    U_1(S) & = \{\repl{w}{a}{d} \mid w \in S, a\in w,w
    \neq a a a ,w \neq a b a\},\\
    U_2(S) & = \{ a x d \mid x \in \mathbf{[d]} \backslash L(S)\},
    \qquad 
    U_3(S) = \{ d x a \mid x \in \mathbf{[d]} \backslash L(S)\},\\
    U_4(S) & = \{ u d v \mid u,v \in L(S)\}, \qquad \;\;\;
    U_5(S) = \{ u d d \mid u\in L(S)\}.
  \end{align*}
\end{definition}


\begin{proposition}
  \label{prop:fismis}
  If $S$ is a maximum independent set of $B(d,3)$, then $f(S)$ is a
  maximum independent set of $B(d+1,3)$.
\end{proposition}
\begin{proof}By definition, $f(S)$ is made up of six disjoint sets. We
  will see that $f(S)$ is an independent set and that it has the right
  cardinality, as in~(\ref{eq:LichBound}).  We start by showing that
  $f(S)$ is an independent set.  This amounts to noticing that there
  are no arrows between the six sets defining $f(S)$.  The only remark
  to bear in mind is that $a x a$ is in $S$ for all $x$, and that $b x
  b$ is also in $S$, except for $x = a$. We leave the details to the
  reader.

  We now compute the cardinality of $f(S)$.  Let $l$ be the number of
  loops of~$S$. We have $|S| = 1 + (d^3-d)/3$, and
  \begin{align*}
   |U_1(S)| &= (d-1)^2 - (l - 1) + (d-l) = d^2 - d +2 - 2l \\
   |U_2(S)| &= |U_3(S)| = (d-l), \qquad |U_4(S)| = l^2, \qquad |U_5(S)| = l.
  \end{align*}
  Only the cardinality of $U_1(S)$ requires explanation. Notice that
  $B(d,3)$ has $(d-1)^2$ cycles whose nodes contain the digit $a$
  once. Each of these contributes one element to $S$ and thus to
  $U_1(S)$, with the exception of $a b b \to b b a \to b a b$ in the
  case that $l = 2$, that contributes no node to $S$ nor
  $U_1(S)$. Likewise, $S$ and $U_1(S)$ contain one element from each
  of the $d-l$ cycles of the form $a a x \to a x a \to x a a$, where
  $x$ is not a loop. Hence, $|U_1(S)|=d^2-d+2-2l$.

  We add the sizes of our six constituents, to obtain 
  \begin{equation*}
    |f(S)| = \frac{(d+1)^3-(d+1)}{3} + 1 + (l-1)(l-2).
  \end{equation*}
  Since $l$ is either $1$ or~$2$ by Proposition~\ref{prop:lessthanthreeloops}, $f(S)$ has the size of an MIS in $B(d+1,3)$.
\end{proof}

We next define another function very similar to $f$ and prove that it
has analogous properties.
\begin{definition}
  \label{def:f'}
  Let $S$ be a maximum independent set of $B(d,3)$. We define
  $f'(S)\subset B(d+1,3)$ as the union of $S$, the sets $U_1(S)$,
  $U_2(S)$, $U_3(S)$, $U_4(S)$ from Definition~\ref{def:f}, and $
  U_5'(S) = \{ d d u \mid u\in L(S)\}$, which is the reverse of $U_5(S)$.
\end{definition}

\begin{proposition}
  \label{prop:f'isgood}
  If $S$ is a maximum independent set of $B(d,3)$, then $f'(S)$ is a
  maximum independent set of $B(d+1,3)$.
\end{proposition}
\begin{proof}
  This proposition is proved analogously to
  Proposition~\ref{prop:fismis}.
\end{proof}

We now define another pair of operators $g$ and $g'$. These will send
a maximum independent set of $B(d,3)$ to a maximum independent set of
$B(d+2,3)$. As before, we follow the convention of Notation~\ref{not:replaceAndLoops}.
\begin{definition}
  \label{def:g}
  Let $S$ be a maximum independent set of $B(d,3)$. We define
  $g(S)\subset B(d+2,3)$
  to be the union  $S \cup \bigcup_{i=1}^7 V_i(S)$, where
  \begin{align*}
    V_1(S) &= \{ \repl{w}{a}{y} \mid y\in \{d,d+1\}, w\in S, a\in w, w\neq
    a a a, w\neq
    a b a\},\\
    V_2(S) &= \{ a x y \mid x\in \mathbf{[d]}\backslash L(S)
    , y \in \{d,d+1 \}\},\\
    V_3(S) &= \{ y x a \mid x\in \mathbf{[d]}\backslash L(S)
    , y \in \{d,d+1 \}\},\\
    V_4(S) &= \{ y x z \mid y,z \in \{d,d+1\}, y\neq z, x\in
    \mathbf{[d]}\backslash L(S)\},\\
    V_5(S) &= \{ u y v \mid u,v \in L(S),
    y \in \{d,d+1\}\},\\
    V_6(S) &= \{ u y y \mid u\in L(S)\}
    , y \in \{d,d+1\}\},\\
    V_7(S) &= \{ y z u \mid y,z \in \{d,d+1\},y\neq z, u\in L(S)\},\\
    V_8(S) &= \{d(d+1)(d+1), (d+1)dd\}.
  \end{align*}
\end{definition}
\begin{definition}
  \label{def:g'}
  Let $S$ be a maximum independent set of $B(d,3)$. We define
  $g'(S)\subset B(d+2,3)$
  to be the union of $S$, the sets $V_1(S)$, $V_2(S)$, $V_3(S)$,
  $V_4(S)$, $V_5(S)$ from Definition~\ref{def:g}, and the sets
  \begin{align*}
    V'_6(S) &= \{ y y u,\ u\in L(S)
    , y \in \{d,d+1\}\},\\
    V'_7(S) &= \{ u y z,\ y,z \in \{d,d+1\},y\neq z, u\in L(S)\},\\
    V'_8(S) &= \{(d+1)(d+1)d, dd(d+1)\},
  \end{align*}
  which are the reverses of $V_6(S)$, $V_7(S)$, and $V_8(S)$ respectively.
\end{definition}

\begin{proposition}
  \label{prop:gg'aremis}
  If $S$ is a maximum independent set of $B(d,3)$, then $g(S)$ and
  $g'(S)$ are
  maximum independent sets of $B(d+2,3)$.
\end{proposition}
\begin{proof} We will prove the statement for the set $g(S)$. The
  result for $g'(S)$ can be proven analogously.  The set $g(S)$ is
  made up of nine disjoint sets. By definition, it is easy to see that
  $g(S)$ is an independent set. We now show that it has the desired
  cardinality. We have $|S| = 1+ (d^3-d)/3$. If $l$ is the number of
  loops of $S$, then
  \begin{align*}
    |V_1(S)| &= 2|U_1(S)| = 2(d^2-d+2-2l),\\
    |V_2(S)| &= 2|U_2(S)| =  2|U_3(S)| = |V_3(S)| = 2(d-l),\\
    |V_4(S)| &= 2|U_4(S)| = 2l^2, \qquad |V_5(S)| = 2(d-l),\\
    |V_6(S)| &= 2|U_5(S)| = 2l, \qquad |V_7(S)| = 2l, \qquad |V_8(S)|
    = 2.
  \end{align*}
  The sum of these sizes is
  $
    |g(S)| = \big((d+2)^3 - (d + 2)\big)/3+ 1 + 2(l-1)(l-2).
    $
  Since $l=1$ or~$2$, the result follows.
\end{proof}


\section{Action of the Symmetric Group on $B(d,3)$}
\label{sec:symmetries-bd-3}

In this section, we study the interaction between $\mathds{S}_d$, the
group of permutations of $\mathbf{[d]}$, and
the four functions we defined in the previous section. In particular,
we show that, up to a permutation of the digits, every maximum
independent set in $B(d,3)$ can be obtained uniquely by successively
composing our four
operators and evaluating this new function at a maximum
independent set of $B(1,3)$ or $B(2,3)$. 

The group $\mathds{S}_d$ acts on the nodes of $B(d,D)$ by $\sigma(x_1
\cdots x_D) = \sigma(x_1)\cdots \sigma(x_D)$ for $\sigma \in
\mathds{S}_d$. This action preserves the graph structure, and
therefore permutes the maximum independent sets.  We will write $A
\sim B$ to mean $A$ and $B$ are two sets in the same orbit under the
action of $\mathds{S}_d$.  Note that the functions $f$, $f'$, $g$,
and~$g'$ are defined so that if $A \sim B$, then $f(A) \sim f(B)$,
etc.  Therefore, each of these functions takes an $\mathds{S}_d$-orbit
of MISs to an $\mathds{S}_{d+1}$- or $\mathds{S}_{d+2}$-orbit of MISs.

\begin{proposition}
  \label{prop:stabf}
  Let $S$ be a maximum independent set of $B(d,3)$. Let $H \subset
  \mathds S_d$ and $H', H'' \subset \mathds S_{d+1}$ 
  be the stabilizers of $S$, $f(S)$ and $f'(S)$,
  respectively. Then $H = H' = H''$,
  where we identify $H$ with its image under the inclusion
  $\mathds{S}_d \hookrightarrow \mathds{S}_{d + 1}$.
\end{proposition}
\begin{proof}
  We only show the equality $H=H'$. The result for $H$ and $H''$ will
  follow  in much the same way.
  We know that \(H \subseteq H'\), and we must prove the other
  inclusion.  Let $\sigma \in H'$, and let $L(S)$ be the loops of $S$,
  with $a$ the distinguished loop with $m_a(S) = 0$. The set of loops
  must be preserved by
  $\sigma$ and moreover, by Proposition~\ref{prop:lessthanthreeloops},
  $\sigma$ fixes each loop.  We want to show that $\sigma(d) = d$.
  Suppose that $\sigma(d) = z \neq d$ and $\sigma(x) = d$, for some
  $x \neq d$. Since $x$ is not a loop, the node $a x d$ then belongs to
  the set $U_2(S)$ from Definition~\ref{def:f}, and so to $f(S)$. That
  means that $\sigma(a x d) = adz$ must be in $f(S)$.  Since this word begins
  with $a$, and has $d$ in the middle, it could only be in~$U_4(S)$. But
  $z \notin L(S)$, and so $adz \notin U_4(S)$. Therefore, $\sigma(d) = d$.

  Now, since $\sigma(d) = d$, $\sigma$ is also an element of
  $\mathds{S}_d$. Furthermore, it must be in the stabilizer of $S$. Otherwise,
  it should map a node of $S$ into a node having a $d$. Since this is
  not possible, $\sigma\in H$.
\end{proof}

\begin{proposition}
  \label{prop:stabg}
  Let $S$ be a maximum independent set of $B(d,3)$. Let $H \subset
  \mathds S_d$ and $H', H'' \subset \mathds S_{d+2}$ be the
  stabilizers of $S$, $g(S)$ and $g'(S)$,
  respectively. Let $\tau \in \mathds S_{d+2}$ be the transposition
  interchanging $d$ and
  $d+1$. Then
  \begin{equation*}
    H' = H'' = \langle\tau,H\rangle,
  \end{equation*}
  where, again, we identify $H$ with its image in $\mathds{S}_{d+2}$.
  Note that $\tau$ commutes with every element of $H$.
\end{proposition}
\begin{proof}Again, we only show the equality $H'=\langle\tau,H\rangle$,
  since the statement involving $H''$ is analogous.  

  As in the proof of Proposition~\ref{prop:stabf}, we know that
  $\langle \tau,H \rangle \subseteq H'.$ Now, let $\sigma \in
  H'$. Again, $\sigma$ must preserve the set $L(S)$ of loops in $g(S)$, and
  by Proposition~\ref{prop:lessthanthreeloops}, $\sigma$ in fact fixes
  each loop. We will show that either $\sigma$ or $\tau\sigma$ fixes
  $d$ and $d+1$. Let $x$, $y$, $z$ and $v$ be such that
  \begin{equation*}
  \xymatrix{
    x\ar@{|->}[r]^{\sigma} & d \ar@{|->}[r]^(0.5){\sigma} & y & \mbox{and} & 
    z\ar@{|->}[r]^(0.4){\sigma} & d+1 \ar@{|->}[r]^(0.6){\sigma} & v.
  }
  \end{equation*}
  We know that $x,y,z,v \notin L(S)$.  Suppose that $x$ is neither $d$
  nor~$d+1$. Then
  we must have $d x a \in V_3(S)$ from
  Definition~\ref{def:g}. The node $\sigma(d x a) = y d a$ 
  has to be in $g(S)$, but it can only be in $V_7(S)$. That means that
  $y = d+1$. Likewise, considering
  \begin{equation*}
    \sigma((d+1)z a) = v(d+1)a,
  \end{equation*} we
  have $v = d$. So $\sigma(d) = d+1$ and $\sigma(d+1) = d$. This
  contradicts our assumption about $x$, and implies that $x = d$ or $d+1$. Analogously, $z = d+1$ or $d$. That means that $\sigma$
  fixes $d$ and $d+1$ or that it transposes them. Therefore, either
  $\sigma$ or $\tau\sigma$ is in~$H$, and so $\sigma \in \langle
  \tau,H \rangle$.
\end{proof}


We now show the precise way in which our functions and $\mathds{S}_d$
interact.


\begin{lemma}
  \label{lem:diffff'difforbits}   \label{lem:diffgg'difforbits}
  Let $S$ and $S'$ be maximum independent sets of $B(d,3)$. Then
    $f(S) \not\sim f'(S')$ and $g(S) \not\sim g'(S')$.
\end{lemma}
\begin{proof}
  We first prove the result for $f$ and $f'$. For contradiction,
  suppose that there is $\sigma\in \mathds{S}_{d+1}$ such that $
  f(S) = \sigma f'(S')$. 
  Let $L(S)$ and $L(S')$ be the loops of $S$ and~$S'$. By
  construction, we have $\sigma L(S') = \sigma L(f'(S')) =
  L(f(S))=L(S)$. Call $a$ and~$a'$ the distinguished loops of $S$
  and~$S'$. By Proposition~\ref{prop:lessthanthreeloops},  we know
  that
  $\sigma(a')=a$.

   Let $x \notin L(S')$ and $y \notin L(S)$ be such that $    \xymatrix{
      x\ar@{|->}[r]^{\sigma} & d \ar@{|->}[r]^(0.5){\sigma} & y.
    }
$ 
  Suppose that $y \neq d$. Then the node $a y d$ is in $U_2(S)$, and
  hence in $f(S)$. Therefore, $\sigma^{-1}(a y d)$ must be in
  $f'(S')$.  But $\sigma^{-1}(a y d) = a' d x$, which cannot be in any
  of the sets that make up $f'(S')$. This implies that $y=d$, hence
  $\sigma(d) = d$.  In other words, $\sigma$ lies in the image of
  $\mathds{S}_d$ in $\mathds{S}_{d+1}$, and so $\sigma f'(S') = f'(\sigma S')$. However,
  $f(S)$ has at least one element of the form $u d d$, and $f'(\sigma
  S')$ has none, so $f(S) \not\sim f'(S')$.

  The proof for $g$ and $g'$ is similar. Namely,
  suppose that there exists $\sigma \in \mathds{S}_{d+2}$ such
  that $g(S) = \sigma g'(S')$. 
 Let $x,z \notin L(S')$, $y,v \notin L(S)$ be such that
  \begin{equation*}
  \xymatrix{
    x\ar@{|->}[r]^{\sigma} & d \ar@{|->}[r]^(0.5){\sigma} & y & \mbox{and} & 
    z\ar@{|->}[r]^(0.4){\sigma} & d+1 \ar@{|->}[r]^(0.6){\sigma} & v.
  }
  \end{equation*}
  Suppose that $y \neq d,d+1$. Then the node $a y d$ is in $V_2(S)$,
  and therefore in $g(S)$. That means that $\sigma^{-1}(a y d) = a' d
  x$ must be in $g'(S')$. But such a node does not belong to any of
  the sets that make up $g'(S')$. This implies that either $\sigma(d)
  = d$ or $\sigma(d) = d+1$.  Analogously, we can prove that
  $\sigma(d+1) = d+1\ \text{ or }\ \sigma(d+1) = d$.

  Therefore, $\sigma$ transposes $d$ and $d+1$ or leaves them
  fixed. By Proposition~\ref{prop:stabg}, the transposition $(d, d+1)$
  is in the stabilizer of $g'(S')$ and so by possibly multiplying
  $\sigma$ on the right by this transposition, we can assume that
  $\sigma$ fixes $d$ and $d+1$ and so it lies in
  $\mathds{S}_d$. Therefore, $\sigma g'(S') = g'(\sigma S')$, but
  $g(S)$ has at least one node of the form $u d d$, and $g'(\sigma
  S')$ has none, so $g(S) \not\sim g'(S')$.
\end{proof}

We now state two invariants that completely characterize maximum
independent sets of $B(d,3)$. This is useful to prove that our
functions $f\!$, $f'\!\!$, $g$, and $g'\!$, together with the action
of $\mathds{S}_d$, allow us to construct all maximum independent sets
of $B(d,3)$. In order to reverse these functions, we make the
following observation, which also holds for loop-less maximum
independent sets. Since we will use it in
Section~\ref{sec:loopless-miss} we state it in full generality.

\begin{proposition}
  \label{prop:intersectmis}
  Let $S$ be a (possibly loop-less) maximum independent set of
  $B(d,3)$, with loops $L(S)$. Let $d'$ be any integer
  such that $c < d' < d$ for all $c \in L(S)$. Then, $    S' = S \cap
  B(d',3)$ 
  is a maximum independent set of $B(d',3)$ with loops $L(S)$.
\end{proposition}
\begin{proof}
  Since $B(d',3)$ is a subgraph of $B(d,3)$, $S'$ is clearly an independent
  set. Furthermore, since $S$ has one element from each cycle except possibly a
  cycle that only uses the digits $a$ and~$b$, then $S'$ has the same property.
  Therefore, $S'$ has the cardinality of a maximum independent set.
\end{proof}

\begin{proposition}
  \label{prop:structff'}
  Let $S$ be a maximum independent set of $B(d,3)$ with $l$ loops,
  where $d$ is at least $3$.
  There exists a digit $x$ such that $m_x(S) = l + 1$
  \emph{if and only if} there exist $\sigma \in \mathds{S}_d$ and $S'$ a
  maximum independent set of $B(d-1,3)$ such that \(S = \sigma f(S')\)
  or \(S = \sigma f'(S')\).
\end{proposition}
\begin{proof}
  The reverse implication follows from the definitions of $f$ and $f'$,
  taking $x = \sigma(d-1)$.
  Conversely, suppose that there is an $x$ with $m_x(S)=l+1$. We know it is not a
  loop by Proposition~\ref{prop:lessthanthreeloops}.
  We define the transposition $\sigma = (d-1,x)$ and the set $S' =
  \sigma S \cap B(d-1,3)$,
  which is a maximum independent set of $B(d-1,3)$ by
  Proposition~\ref{prop:intersectmis}.

  Let $a$ denote the distinguished loop of $S$. We know that the node $x a x
  \notin S$. Therefore, either $x x a$ or $a x x$ must be in
  $S$.
  Suppose that $a x x \in S$, in which case we claim that $S = \sigma f(S')$.

  We now consider each of the sets that make up $\sigma
  f(S')$, and show that they are included in~$S$.
  The nodes of $\sigma S'$ belong to $S$, by definition of $S'$.
  Let us consider the nodes of $\sigma U_1(S')$. The
  nodes of this set are of the form $x y x$, $x y y$, $y y x$, $x y z$
  or $y z x$, for $y$ and~$z$ distinct from~$x$ and $y, z \notin
  L(S)$.
  \begin{itemize}
  \item The nodes of the form $x y x$ are all in $S$ by the hypothesis
    on~$x$.
  \item If $x y y \in \sigma U_1(S')$, then $a y y \in S'$. This means
    that $a y y \in S$, and so $y y x$ cannot be in $S$. The node $y x y$
    cannot be in $S$ either, since $x y x$ is. So, $x y y \in
    S$. Analogously, if $y y x \in \sigma U_1(S')$, then $y y x \in S$.
  \item If $x y z \in \sigma U_1(S')$, then $a y z \in S'$ and $a y z \in
    S$. Since neither $z x y$ (adjacent to $x y x$) nor $y z x$ (adjacent to
    $a y z$) can be in $S$, $x y z$ must be in $S$. The same reasoning
    applies to $y z x$.
  \end{itemize}

  Let us consider the nodes of $\sigma U_2(S')$. These have the
  form $a y x$. The nodes $y x a$ (adjacent to $x y x$) and $x a y$ (adjacent
  to $a y a$) cannot be in $S$, which implies that $a y x \in S$. The same
  reasoning shows that $\sigma U_3(S') \subset S$.

  A node from $\sigma U_4(S')$ is of the form
  $u x v$, with $u$ and~$v$ loops. The nodes $x u v$ (adjacent to $u x u$) and
  $u v x$ (adjacent to $v x v$) cannot be in $S$. Therefore, $u x v \in S$,
  and $\sigma U_4(S') \subset S$.

  Finally, we know that $a x x \in S$ or $xxa \in S$.  Assume the
  first case. If $S$ has a single loop, we have that 
  $\sigma U_5(S')\subset S$. If $S$ has an
  extra loop $b$, the nodes $x b x$ (adjacent to
  $b x b$) and $x x b$ (adjacent to $a x x$) cannot be in $S$. That implies
  that $b x x \in S$, which means $\sigma U_5(S') \subset S$.
  This proves that $S \supseteq \sigma f(S')$. Since both sets have
  the same cardinality, we 
  conclude that equality holds.

  On the other hand, if $x x a \in S$, an analogous
  procedure shows that $S = \sigma f'(S')$.
\end{proof}

The following lemma is used in the proof of
Proposition~\ref{prop:structgg'},
which is the analogue of
Proposition~\ref{prop:structff'} for the operators $g$ and $g'$.
Note that in preparation for our study of loop-less maximum independent
sets in Section~\ref{sec:loopless-miss}, we prove
Lemma~\ref{lem:twobigmaxima} for loop-less maximum independent sets as
well.
\begin{lemma}
  \label{lem:twobigmaxima}
  Let $S$ be a (possibly loop-less) maximum independent set of
  $B(d,3)$, with $d \geq 3$. If there exist two different digits $y$
  and $z$, which are not loops, such that
  \begin{equation*}
    m_y(S) = m_z(S) = l + 2,
  \end{equation*}
  then $y z y \notin S$ and $z y z \notin S$.
\end{lemma}
\begin{proof}
  Suppose that $y z y \in S$. Then, by the assumptions on $m_y(S)$,
  there must be some $v \neq y$ such that $y v y \notin S$. Suppose
  that $v y y \in S$. The node $z y z$ cannot be in $S$, and by the
  assumption on $m_z(S)$, $z v z \in S$. Therefore, the nodes $z v y$
  (adjacent to $v y y$), $v y z$ (adjacent to $y z y$) and $y z v$
  (adjacent to $z v z$) are not in $S$. But then the cycle $z v y \to
  v y z \to y z v$ contributes no nodes to $S$, which contradicts the
  fact that $S$ has maximum cardinality. If we assume that $y y v
  \in S$, then the cycle $y v z \to v z y \to z y v$ cannot contribute
  any node to~$S$, a contradiction.

  In conclusion, our assumption that $y z y$ is in $S$ is inconsistent
  with $S$ being a maximum independent set. By symmetry, the same
  holds if we assume $z y z \in S$.
\end{proof}

\begin{proposition}
  \label{prop:structgg'}
  Let $S$ be a maximum independent set of $B(d,3)$, $d\geq 3$, with
  $l$ loops ($l=1$ or $2$). Then, there are two different digits $y$
  and $z$ such that
  \begin{equation*}
    m_y(S) = m_z(S) = l + 2
  \end{equation*}
  and \emph{no} digit $x$ such that \(m_x(S) = l + 1\), \emph{if
    and only if} there exist $\sigma \in \mathds{S}_d$ and $S'$ a maximum
  independent set of $B(d-2,3)$ such that
  \begin{equation*}
    S = \sigma g(S')\ \text{ or }\ S = \sigma g'(S').
  \end{equation*}
\end{proposition}
\begin{proof}
  One implication follows from the construction of $g$ and $g'$ taking
  $y = \sigma(d-1)$ and $z = \sigma(d-2)$.  The proof in the other
  direction is analogous to the proof of
  Proposition~\ref{prop:structff'}. We can safely assume that $y =
  d-1$ and $z = d-2$. By Lemma~\ref{lem:twobigmaxima}, either
  the pair $(d-1)(d-2)(d-2)$ and $(d-2)(d-1)(d-1)$ are in~$S$, or
  the pair $(d-1)(d-1)(d-2)$ and $(d-2)(d-2)(d-1)$ are in~$S$. In the former
  case, we find that there is an $S'$ such that $S = \sigma g(S')$. In
  the latter case, we find that $S = \sigma g'(S')$.
\end{proof}

\begin{corollary}
  \label{cor:fandgdonooverlap}
  Let $S$ and $S'$ be maximum independent sets of $B(d-1,3)$ and
  $B(d-2,3)$ with $d \geq 3$. Then for $\mathcal{F} = f,f'$ and
  $\mathcal{G} = g,g'$, we have $    \mathcal{F}(S) \not \sim \mathcal{G}(S')$.
\end{corollary}
\begin{proof}
  This result follows from the invariants of $\mathcal{F}(S)$ and
  $\mathcal{G}(S')$ that are stated in
  Propositions~\ref{prop:structff'} and~\ref{prop:structgg'}.
\end{proof}

This corollary, together with Lemmas~\ref{lem:diffff'difforbits}
and~\ref{lem:diffgg'difforbits}, shows that all four functions produce
essentially different (i.e.\ in different $\mathds{S}_d$-orbits)
maximum independent sets.

\section{Characterization of Maximum Independent Sets}
\label{sec:maintheorem}

In this section, we show that the functions $f$, $f'$, $g$, and $g'$,
together with the action of $\mathds{S}_d$ are sufficient to construct
every maximum independent set of $B(d, 3)$.  For the rest of this
section, $L$ will denote the set of loops of $S$, and $l$ will denote
the cardinality of~$L$. In Section~\ref{sec:loopless-miss}, we will work with loop-less
maximum independent sets. For that reason, we prove some of the
results of this section in that context too.

As we mentioned in Section~\ref{sec:second-results}, the sets
$\M_x(S)$ from Definition~\ref{def:M_x(A)} play a key role. We start
our discussion with two technical lemmas about them.

\begin{lemma}
  \label{lem:threemaxima} \label{lem:twoandonemaxima} Let $S$ be a
  (possibly loop-less) maximum independent set of $B(d,3)$. There
  cannot be three different digits $x$, $y$, and $z$, with $x,y,z
  \notin L$, 
  such that
  \begin{equation}
    \label{eq:3}
    \begin{aligned}
    \M_x(S) &= \M_y(S) = L \cup \{x,y,z\},\\
    \M_z(S) &= L \cup \{x,z\} \mbox{ or } L \cup \{x,y,z\}.
    \end{aligned}
  \end{equation}
\end{lemma}
\begin{proof}
  Suppose that $S$ is a maximum independent set and $x$, $y$, and $z$
  satisfy~\eqref{eq:3}.  Without loss of generality, we can assume
  that $x$, $y$, $z$, and the loops are smaller than $l + 3$. Then $S' =
  S \cap B(l+3, 3)$ is a maximum independent set in $B(l+3,3)$ by
  Proposition~\ref{prop:intersectmis} with $\M_x(S') = \M_x(S)$,
  $\M_y(S') = \M_y(S)$, and $\M_z(S') = \M_z(S)$.

  Without loss of generality we may assume that $x y y,\allowbreak y x
  x,\allowbreak x z z,\allowbreak z y y,\allowbreak z x x \in S'$
  since $y x y,\allowbreak x y x,\allowbreak z x z,\allowbreak y z
  y,\allowbreak x z x \notin S'$. But this implies that there is no
  element of the cycle containing $z y x$ in $S'$, a contradiction.
  Therefore,  no such $S$ exists.
\end{proof}

\begin{lemma}
  \label{lem:threemaximabis}
  Let $S$ be a (possibly loop-less) maximum independent set of
  $B(d,3)$. There cannot be three different digits $x$, $y$, and $z$,
  none of which are loops, such that
  \begin{equation*}
    m_x(S) = m_y(S) = m_z(S) = l + 2.
  \end{equation*}
\end{lemma}
\begin{proof}
  We prove the result by contradiction. Suppose there are such $x$,
  $y$ and $z$. We know that $L \cup \{x\} \subset \M_x(S)$ and
  $|\M_x(S)| = l + 2$. Therefore, at least one of $y$ and $z$ must
  appear between~$x$. An analogous statement holds for $y$ and
  $z$. Without loss of generality, suppose that $y$ appears
  between~$x$. Then $y x y$ (adjacent to $x y x$) is not in $S$, which
  forces
  $z$ to appear between $y$. That, in turn, forces $x$ to appear
  between~$z$. That is, the nodes $x y x$, $y z y$ and $z x z$ are in
  $S$. But then, none of the nodes $x y z \to y z x
  \to z x y$ are in $S$, contradicting the
  maximality of $S$.
\end{proof}

\begin{remark}
  \label{rem:twomaxima} Note that a maximum independent set $S$ of
  $B(d,3)$ with $l$ loops can have at most one digit satisfying
  $m_x(S)=l+1$. If there were two, say $x$ and $y$, then $x y x$ and
  $y x y$ would have to be in $S$, a contradiction.
\end{remark}

The next proposition shows
that, up to permutation, any maximum independent set lies in the
image of one of our four operations.
\begin{proposition}
  \label{prop:mainlemma}
  Let $S$ be a (possibly loop-less) maximum independent set of
  $B(d,3)$ with $d \geq 3$. Suppose there is no digit $z$ such that
  $m_z(S) = l + 1$.  Then, there must be exactly two digits $x$ and
  $y$ such that $m_x(S) = m_y(S) = l + 2$. Moreover, $\M_x(S) =
  \M_y(S) = L \cup \{x,y\}$.
\end{proposition}
\begin{proof}
  We just need to show that $m_x(S) = m_y(S) = l + 2$.
  Lemma~\ref{lem:twobigmaxima} 
  implies that $\M_x(S) = \M_y(S) = L \cup \{x,y\}$.
  By reordering the digits, we can assume that $m_{d-1}(S) \leq
  m_{d-2}(S)
  \leq m_i(S)$ for all $i < d-2$.  By hypothesis, we know that
  $m_{d-1}(S) \geq l + 2$ and we want to prove that $m_{d-2}(S) = l +
  2$. Lemma~\ref{lem:threemaximabis} will then imply that $d-2$
  and $d-1$ are the only digits with this property.

  We prove that $m_{d-2}(S)=l+2$ by induction on $d$. Our base cases are $d\leq
  l+3$. If $d = l + 1$, then the unique $z$ not in~$L$ satisfies $m_z(S) = l +
  1$, which contradicts our hypothesis. If $d = l + 2$, and $x$
  and~$y$ are not in~$L$, then $m_x(S) \geq l + 2$ implies $m_x(S) = l
  + 2$, and likewise for $y$. If $d = l + 3$, then
  Lemma~\ref{lem:twoandonemaxima} gives us the result.

  Now, let $d$ be greater than $l + 3$ and consider $S' = S \cap
  B(d-1,3)$. By the inductive hypothesis, we must have one of two
  possibilities:
  
  \textbf{Case 1:} $S'$ has exactly one digit $z$ with $m_z(S') =
  l+1$.  If $z = d-2$, we are done. Suppose that $z \neq d-2$. By
  Remark~\ref{rem:twomaxima}, $m_{d-2}(S') > m_z(S')$, $m_{d-2}(S)
  \leq m_z(S)$ and  $m_{d-2}(S')\leq 
  m_{d-2}(S)$. Thus, we must
  have $m_z(S') = m_z(S)-1$ and $m_{d-2}(S') = m_{d-2}(S)$. This means
  that $z(d-1)z$ is not in~$S$ and $(d-2)(d-1)(d-2) \in S$, which  implies 
  that $m_{d-2}(S) = l + 2$, as we wanted to show.
  
  \textbf{Case 2:} $S'$ has exactly two digits $x$ and $y$ with $m_x(S') =
  m_y(S') = l + 2$. We split this situation in two subcases.
  
  \textbf{Case 2.1:} We suppose $x,y \neq d-2$. By an argument similar
  to that of Case 1, we know that $\M_x(S) =
  \M_y(S) = L \cup \{x,y,d-1\}$ and
  \[
  m_{d-2}(S) = l + 3, \qquad \M_{d-2}(S)
  \supseteq L \cup \{d-2, d-1, x, y\},
  \]
  which is a contradiction.


  \textbf{Case 2.2:} Either $x$ or $y$ equals $d-2$. Suppose $y =
  d-2$. Since $m_{d-2}(S') = l + 2$, then $m_{d-2}(S) = l + 2$ (and
  we are done) or
  $m_{d-2}(S) = m_x(S) = l + 3$. Hence,
%
  \begin{equation}
    \label{eq:0}
    \M_x(S) = \M_{d-2}(S) = L \cup \{x,d-2,d-1\}.
  \end{equation}
  Since $m_{d-1}(S) \leq m_{d-2}(S) = l+3$, we have that
  \[
  \M_{d-1}(S) = L \cup \{d-1, u, v\} \text{ or } L \cup \{d-1,u \}.
  \]

  \textbf{Case 2.2.1:} Suppose $m_{d-1}(S) = l+3$. We will show that
  $\M_{d-1}(S) = L \cup \{x,d-2,d-1\}$, which, together
  with~\eqref{eq:0}, contradicts Lemma~\ref{lem:twoandonemaxima}.

  Assume $u,v \neq d-2$. That means that $(d-1)(d-2)(d-1) \in
  S$. Since $u\neq v$, we can assume without loss of generality that
  $u \neq x$. Then $x u x \in S$ and $(d-2)u(d-2) \in S$. The nodes
  $(d-1)u(d-2)$ and $(d-2)u(d-1)$ must be in $S$, because the rest of
  the nodes in their cycles are adjacent to something just shown to be
  in $S$. We know that $(d-1)u(d-1)$ is not in~$S$, because of the
  definition of~$u$. Additionally, the nodes $(d-1)(d-1)u$ and $u(d-1)(d-1)$
  are adjacent to the nodes we just showed are in~$S$. Therefore,
  neither of them belong to $S$, a contradiction.
Hence, one of $u$ and~$v$ must equal $d-2$, and so we have
  \begin{equation*}
    \M_{d-1}(S) = L \cup \{u,d-2,d-1\}.
  \end{equation*}

  To finish, we need to prove that $u = x$. Assume the contrary. Then
  $x u x$ and $(d-1)x(d-1)$ are in $S$. Therefore, by inspecting their
  cycles we see that both $x u (d-1)$ and $(d-1)u x$ must be in~$S$.
  On the other hand, either $u(d-1)(d-1) \in S$ or $(d-1)(d-1)u \in
  S$.  However, $u(d-1)(d-1) \in S$ implies $x u (d-1)\notin S$, and
  $(d-1)(d-1)u \in S$ implies $(d-1) u x \notin S$. Therefore, $u =
  x$.

  \textbf{Case 2.2.2} Suppose $m_{d-1}(S) = l+2$. If we assume $x$ and
  $d-2$ are not in $\M_{d-1}(S)$ and proceed as in the previous case,
  we get a contradiction. Therefore, Lemma~\ref{lem:twoandonemaxima}
  applied to $x$, $d-1$ and $d-2$ leads to a contradiction.
\end{proof}

We now state our main result.

\begin{theorem}[\textbf{Characterization of the Maximum Independent
    Sets of $\mathbf{B(d,3)}$}]
  \label{theo:main}
  For all positive $d$ we have:
  \begin{enumerate}
  \item \label{item:unique-sequence} Any orbit of independent sets of
    $B(d,3)$ under the action of $\mathds{S}_d$ is obtained from the
     $\{0 0 0\}$ and the orbit of $\{0 0
    0 ,0 1 0,1 1 1\}$ under $\mathds{S}_2$ by a \emph{unique} sequence
    of applications of $f\!$, $f'\!\!$, $g$, and $g'\!$.
  \item \label{item:transp} Let $S$ be an MIS of $B(d,3)$. Then the
    subgroup of $\mathds S_d$ stabilizing $S$ is generated by disjoint
    transpositions. In particular, the cardinality of the
    stabilizer of~$S$ is a power of~$2$.
  \item \label{item:b_dk} Let $b_{d,k}$ be the number of orbits of
    MISs in $B(d,3)$ whose elements have stabilizers of size
    $2^k$. Then we have the recurrence relation
    \begin{equation*}
      \begin{cases}
        b_{1,0} = 1, \;
        b_{2,0} = 3,\\
        b_{d,k} = 2b_{d-1,k} + 2b_{d-2,k-1} & \text{for } d \geq 3,
      \end{cases}
    \end{equation*}
    and the generating function
    \begin{equation*}
      \sum_{d=1}^\infty \sum_{k=0}^\infty b_{d,k} t^d s^k =
      \frac{t+t^2}{1-2t-2t^2s}.
    \end{equation*}
  \item \label{item:a_d} The number $a_d$ of maximum
  independent sets of $B(d,3)$ satisfies 
    \begin{equation*}
      \begin{cases}
        a_1 = 1,\;
        a_2 = 6,\\
        a_d = 2d a_{d-1} + d(d-1)a_{d-2} & \text{for } d \geq 3,
      \end{cases}
    \end{equation*}
  and has exponential generating function
  \begin{equation*}
    \sum_{d=1}^\infty \frac{a_d t^d}{d!} = \frac{t+t^2}{1-2t-t^2}.
  \end{equation*}
  \end{enumerate}
\end{theorem}
\begin{proof}
  For $d=1$, the only maximum independent set of $B(1,3)$ consists of
  the unique node $\{000\}$. For the case of $d=2$, it can be checked
  manually that the three orbits of maximum independent sets under
  $\mathds{S}_2$ are the orbits of $\{000, 010, 011\}$, $\{000, 010,
  110\}$, and $\{000, 010, 111\}$. Note that the first two of these
  are $f(\{000\})$ and $f'(\{000\})$ respectively.  Thus, the
  existence statement in (\ref{item:unique-sequence}) follows from
  Propositions~\ref{prop:structff'}, \ref{prop:structgg'}
  and~\ref{prop:mainlemma}. The uniqueness comes from
  Lemma~\ref{lem:diffff'difforbits} and
  Corollary~\ref{cor:fandgdonooverlap}.

  The statements in (\ref{item:transp}) and~(\ref{item:b_dk}) follow
  from the previous result and the description of the stabilizers in
  Propositions~\ref{prop:stabf} and~\ref{prop:stabg}. Finally, the
  generating function in (\ref{item:a_d}) is obtained by substituting
  $s=1/2$ into the previous generating function, because
  \begin{equation*}
  a_d = \sum_{k=0}^\infty \frac{d! b_{d,k}}{2^k}.
  \end{equation*}
  The recurrence follows immediately.
\end{proof}

The following table lists the values of $b_{d,k}$, for all $d \leq 6$.

\vspace{1ex}
\begin{center}
  \begin{tabular}{|c|c|c|c|c|c|c|c|}
    \hline
    $k \backslash d$ & 1 & 2 & 3 & 4 & 5 & 6\\
    \hline
    \multirow{2}{*}{0} & 1 & 3 & 6 & 12 & 24 & 48\\
    & (1,0) & (2,1) & (4,2) & (8,4) & (16,8) & (32,16)\\
    \hline
    \multirow{2}{*}{1} & & & 2 & 10 & 32 & 88\\
    & & & (2,0) & (8,2) & (24,8) & (64,24)\\
    \hline
    \multirow{2}{*}{2} & &  &  & & 4 & 28\\
    & & & & & (4,0) & (24,4)\\
    \hline
  \end{tabular}
\end{center}
\vspace{1ex}

\noindent In each entry, the first number indicates the number of orbits whose
elements have only one loop. The second one is the number of orbits
with two loops.

\section{Loop-less Maximum Independent Sets}
\label{sec:loopless-miss}

In this section, we analyze the number of \emph{loop-less} maximum
independent sets (LMISs) of $B(d,3)$, for all $d$. Recall from the
introduction that the size of an LMIS of $B(d,3)$ is
\begin{equation*}
  \alpha^*(d,3) = \frac{d^3 - d}{3} = \alpha(d,3) - 1.  
\end{equation*}
By MIS, we will continue to mean a maximum independent set \emph{with}
loops. As in previous sections, we let $a$ be the loop of $S$ such
that $m_a(S)= 0$, and the other loop (if there is one) is denoted by
$b$.

In what follows, we provide an explicit bijection between LMISs and
MISs of $B(d,3)$. 

\begin{definition}
  Let $S$ be a maximum independent set of $B(d,3)$, $d \geq 3$. We
  define
  \begin{equation}
    \label{eq:6}
    h(S) =
    \begin{cases}
      S \backslash \{a a a\} & \text{if $S$ has only one loop,}\\
      S \backslash \{a a a,b b b,a b a\} \cup \{a a b,b b a\} &
      \text{if $S$ has two loops $a < b$,}\\
      S \backslash \{a a a,b b b,a b a\} \cup \{b a a,a b b\} &
      \text{if $S$ has two loops $a > b$.}
    \end{cases}
  \end{equation}
\end{definition}

\begin{proposition}
  \label{prop:hgood}
  Let $S$ be a maximum independent set of $B(d,3)$. Then $h(S)$ is an
  LMIS of $B(d,3)$.
\end{proposition}
\begin{proof}
  Let $S$ be an MIS of $B(d,3)$. If $S$ has only one loop, then
  eliminating it leaves us with an independent set of the correct
  size.

  If $S$ has two loops, say $a$ and $b$, then $h(S)$ is a set of the
  correct size, since the nodes we added were not already present in
  $S$. However, we must see that $h(S)$ is an independent set. Assume
  $a < b$. 
Suppose we have a node adjacent to~$a a b$. Then it is of
  the form $a b x$ or $x a a$. Since $b x b$ and $a a a$ are in $S$,
  then $a b x$ and $x a a$ cannot be in $S$.
A similar argument show that adding $bba$ to $S$
 preserves independence. 
Therefore,
  the nodes we add are not adjacent to any other nodes in the
  construction, and the result follows. The case $a > b$ is proved
  analogously.
\end{proof}

\begin{proposition}
  \label{prop:hinjective}
  The function $h$ is injective.
\end{proposition}
\begin{proof}
  Let $S$ and $S'$ be two different MISs of $B(d,3)$. Then showing
  that $h(S) \neq h(S')$ is just a matter of analyzing all the
  possible combinations of loops and their relative order in $S$ and
  $S'$. We leave the details to the reader.
%
%
%
%
%
%
\end{proof}


\begin{lemma}
  \label{lem:twoloopsinvariant}
  Let $S$ be a maximum independent set with two loops $a$ and $b$. Let
  $\tau$ be the transposition of $a$ and $b$. Let $S' = S \backslash
  \{a a a,b b b,a b a\}$. Then $S' = \tau S'$.
\end{lemma}
\begin{proof}
  We must show that for every node $w \in S'$ such that $a \in w$, we
  have $\repl{w}{a}{b} \in S'$ and vice versa.  Notice that any node
  of $S'$ cannot contain $a$ and~$b$ simultaneously. The nodes that
  contain two $a$'s or two $b$'s are $a x a$ and $b x b$, and they are
  in $S'$ for all $x \neq a,b$. Thus, $x a y \notin S'$ for all $x, y
  \neq a$. 

The nodes that contain only one $a$ are $x y a$ or $a x
  y$ for $x,y \neq a$. If $x y a \in S'$, then $b x y \notin S'$, and
  so $x y b$ must be in $S'$ in order to have one element from its
  cycle. We can prove that $a x y \in S'$ implies $b x y \in S'$ in a
  similar way.
\end{proof}

\begin{proposition}
  \label{prop:surjective}
  The function $h$ is surjective.
\end{proposition}
\begin{proof}
  Let $S$ be an LMIS of $B(d,3)$. By Proposition~\ref{prop:mainlemma},
  we have two possibilities: 

  First, if there is a digit $x$ such that $m_x(S) = 1$, then there is no
  node of the form $x x y$ or $y x x$. Therefore, $S' = S\cup \{x x
  x\}$ is an MIS of $B(d,3)$, and $S = h(S')$.

  Second, if there are two digits $x$ and $y$ such that
  $m_x(S) = m_y(S) = 2$, then we have either $x x y, y y x \in S$ or
  $y x x, x y y \in S$. In the first case, we construct
  \begin{equation*}
    S' = S \cup \{x x x,y y y,x y x\} \backslash \{x x y, y y x\}.
  \end{equation*}
  If $x < y$, then $S = h(S')$. If $x > y$, then by
  Lemma~\ref{lem:twoloopsinvariant}, $S = h(\tau S')$, where $\tau$ is
  the transposition of $x$ and $y$. The remaining case is dealt with
  analogously.
\end{proof}

\begin{theorem}
  \label{theo:loop-less}
  Let $a^*_d$ be the number of loop-less maximum independent sets of
  $B(d,3)$. Then $a^*_d = a_d$.
\end{theorem}
\begin{proof}
  This follows from Propositions~\ref{prop:hinjective} and~\ref{prop:surjective}.
\end{proof}
We conclude with a result that links comma-free codes and loop-less
maximum independent sets.
\begin{theorem}
  \label{theo:comma-free}
  Every loop-less maximum independent set is a maximum comma-free code of
  length~$3$. In particular, the number of equivalence classes of comma-free
  codes in an alphabet of size~$d$ is at least $2^d$, where
  equivalence means equivalence under the action of $\mathds S_d$.
\end{theorem}

\begin{proof}
  If $S$ is an LMIS with $x_1 x_2 x_3$ and $y_1 y_2 y_3$ elements of~$S$, then
  $x_2 x_3 y_1$ cannot be an element of~$S$ because it is adjacent to $x_1 x_2
  x_3$ in $B(d, 3)$. Likewise, $x_3 y_1 y_2$ cannot be in~$S$ because it is adjacent to
  $y_1 y_2 y_3$. Therefore, $S$ is a comma-free code.

  For the second statement, by considering only the first term of the recurrence
  relation in Theorem~\ref{theo:main}(\ref{item:a_d}), we see that $a_d \geq
  2^d d!$. Therefore, the number of maximum comma-free codes is at least $2^d
  d!$, so the number of equivalence classes under the action of~$\mathds S_d$
  must be at least $2^d$.
\end{proof}

The set $S = \{100, 110\}$ is an example of a maximum comma-free code which is
not an independent set for $d = 2$.

\section*{Acknowledgements}

We would like to thank Ver\'onica Becher for introducing de Bruijn
graphs to us. We also thank Alicia Dickenstein for her useful
comments. We are indebted to Yukio Shibata, for kindly sharing the
work of his group with us.
We used the free computer algebra system \cocoa~\cite{cocoa} for the
initial computations which suggested our results.

\bibliographystyle{amsplain}
\bibliography{journals,debruijn}

\end{document}